\documentclass[12pt]{amsart}

\newtheorem{theorem}{Theorem}[section]
\newtheorem{proposition}[theorem]{Proposition}
\newtheorem{lemma}[theorem]{Lemma}
\newtheorem{cor}[theorem]{Corollary}
\newtheorem*{theo}{Theorem \ref{bound}}
\newtheorem*{coro}{Corollary \ref{conjecture}}
\newtheorem*{prop}{Proposition \ref{prop}}

\theoremstyle{definition}
\newtheorem{definition}[theorem]{Definition}

\newtheorem*{rem}{[Remark]}

\newcommand{\Gin}{\ensuremath{\mathrm{Gin}}}

\newcommand{\init}{\ensuremath{\mathrm{in}}}
\newcommand{\GL}{\ensuremath{GL_n (K)}}
\newcommand{\grev}{Gr\"obner basis}
\newcommand{\rank}{\ensuremath{\mathrm{rank}}}
\newcommand{\lex}{{\prec_{rev}}}
\newcommand{\revlex}{{\prec_{rev}}}
\newcommand{\elex}{ \preceq_{rev}}
\newcommand{\rlex}{ \succeq_{rev}}
\newcommand{\sub}{\ensuremath{M_{\succeq S}(J,d)} }
\newcommand{\tsub}{\ensuremath{M_{\succ S}(J,d)} }
\newcommand{\lsub}[1]{\ensuremath{ \widetilde M_{\succeq S} (J,d,{#1})}}

\author{Satoshi Murai}
\address{
Department of Pure and Applied Mathematics\\
Graduate School of Information Science and Technology\\
Osaka University\\
Toyonaka, Osaka, 560-0043, Japan\\}
\email{s-murai@ist.osaka-u.ac.jp}
\thanks{The author is supported by JSPS Research Fellowships for Young Scientists}

\begin{document}
\title[generic initial ideals and  exterior shifting of join]{Generic initial ideals and exterior algebraic shifting of the join  of simplicial complexes}
\maketitle

\begin{abstract}
In this paper,
the relation between algebraic shifting and join which was conjectured by Eran Nevo will be proved.
Let $\sigma$ and $\tau$ be simplicial complexes
and $\sigma * \tau$ their join.
Let $J_\sigma$ be the exterior face ideal of $\sigma$
and $\Delta(\sigma)$ the exterior algebraic shifted complex of $\sigma$.
Assume that  $\sigma * \tau$ is a simplicial complex on $[n]=\{1,2,\dots,n\}$.
For any $d$-subset $S \subset [n]$,
let $m_{\elex S}(\sigma)$ denote the number of $d$-subsets $R\in \sigma$
which is equal to or smaller than $S$ w.r.t.\ the reverse lexicographic order. 
We will prove  that
$m_{\elex S}(\Delta({\sigma * \tau}))\geq m_{\elex S}(\Delta ( {\Delta( \sigma)} * {\Delta (\tau)} ))$
for all $S \subset [n]$.
To prove this fact,
we also prove that
$m_{\elex S}(\Delta(\sigma))\geq m_{\elex S}(\Delta ( {\Delta_{\varphi}( \sigma)} ))$
for all $S\subset [n]$ and for all non-singular matrices $\varphi$,
where $\Delta_{\varphi}(\sigma)$ is the simplicial complex defined by $J_{\Delta_{\varphi}(\sigma)}=\init(\varphi(J_\sigma))$.
\end{abstract}

\section*{Introduction}
Algebraic shifting,
which was introduced by Kalai,
is a map which associates with each simplicial complex $\sigma$
another simplicial complex $\Delta(\sigma)$ with special conditions.
%
%
Nevo \cite{N} studied some properties of algebraic shifting with respect to basic constructions of simplicial complexes,
such as union, cone and join.
With respect to union and cone,
algebraic shifting behaves nicely.
However, with respect to join,
Nevo found that algebraic shifting does not behave nicely contrary to a conjecture by Kalai \cite{K}.

First, we will recall Kalai's conjecture and Nevo's counter example.
Let $\sigma$ be a simplicial complex on $\{1,2,\dots,k\}$,
$\tau$  a simplicial complex on $\{k+1,k+2,\dots,n\}$,
and let $\sigma * \tau$ denote their join,
in other words,
$$\sigma * \tau = \{ S\cup R : S\in \sigma \ \mbox{and} \ R\in \tau\}.$$
Kalai  conjectured that $\Delta(\sigma * \tau)=\Delta(\Delta(\sigma) * \Delta (\tau))$,
where $\Delta (\sigma)$ is the exterior algebraic shifted complex of $\sigma$.
However, Nevo \cite{N} found a counter example.
We quote his example.
Let $\Sigma (\sigma)$ denote the suspension of $\sigma$,
i.e.,
the join of $\sigma$ with two points.
Nevo showed that
if $\sigma$ is the simplicial complex generated by $\{1,2\}$ and $\{3,4\}$ then
the $2$-skeleton of $\Delta(\Sigma(\sigma))$ is $\{\{1,2,3\},\{1,2,4\},\{1,2,5\},\{1,2,6\}\}$
and that of $\Delta(\Sigma(\Delta(\sigma)))$ is $\{\{1,2,3\},\{1,2,4\},\{1,2,5\},\{1,3,4\}\}$.

Next, we will recall Nevo's conjecture.
Let $\sigma$ and $\tau$ be simplicial complexes
on  $[n]=\{1,2,\dots,n\}$ and $\lex$ the reverse lexicographic order $\revlex$ induced by $1<2<\cdots<n$.
In other words,
for  $S\subset [n]$ and $R\subset [n]$ with $S\ne R$,
define $S\lex R$ if 
(i) $|S|<|R|$ or
(ii) $|S|=|R|$ and the minimal integer in the symmetric difference 
$(S\setminus R )\cup (R\setminus S)$
belongs to $S$.
For an integer $d \geq 0$,
we write $\sigma_d=\{ S\in \sigma : |S|=d+1\}$.
Define $\sigma \leq_{Rev} \tau$ if
the smallest element  w.r.t.\ $\revlex$
in the symmetric difference 
between $\sigma_d$ and $\tau_d$ belongs to $\sigma_d$ for all $d \geq 0$,
i.e.,
$\min_{\lex}\{S: S \in ( \sigma_{d} \setminus \tau_d) \cup ( \tau_d \setminus \sigma_d)\} \in \sigma$
for all $d \geq 0$.

Nevo conjectured that (\cite[Conjecture 6.1]{N}), for any simplicial complex $\sigma$, one has
$$\Delta(\Sigma (\sigma)) \leq_{Rev} \Delta (\Sigma (\Delta(\sigma))).$$
(In the previous example,  the symmetric difference  is $\{ \{1,2,6\}, \{1,3,4\}\}$
and $\{1,2,6\} \in \Delta(\Sigma (\sigma))$.)
In this paper, we will prove a stronger result.
For any subset $S \subset [n]$, let 
$$m_{\elex S}(\sigma)=|\{R\in \sigma : \ |R|=|S|\ \mbox{and}\ R \elex S\}|.$$
We will prove the following.
(The definition of $\Delta_{\varphi}(\sigma)$ will be given in \S 3.)

\begin{theo}
Let $\sigma$ be a simplicial complex on $[n]$ and
$\varphi \in \GL$.
Then, for any $S\subset [n]$,
one has
$$m_{\elex S}(\Delta(\sigma)) \geq m_{\elex S}(\Delta (\Delta_{\varphi}(\sigma))).$$
\end{theo}

By using Theorem \ref{bound},
we can easily prove the next corollary
which implies Nevo's conjecture.
\begin{coro}
Let $\sigma$ be a simplicial complex on $\{1,2,\dots,k\}$ and
$\tau$ a simplicial complex on $\{k+1,k+2,\dots,n\}$.
Then, for any $S \subset [n]$, one has 
$$m_{\elex S} (\Delta(\sigma * \tau))\geq m_{\elex S} (\Delta(\Delta(\sigma) * \Delta(\tau))).$$
\end{coro}

Now, we will explain why Corollary \ref{conjecture} implies Nevo's conjecture.
Let $d\geq 0$ be a positive integer,
and let $\mathcal{L} \subset { [n] \choose d}$ and $\mathcal{R} \subset { [n] \choose d}$
be families of $d$-subsets of $[n]$
which satisfy $m_{\elex S}( \mathcal{L}) \geq m_{\elex S} (\mathcal{R})$ for all $S \in {[n] \choose d}$.
Set $T= \min_{\revlex}\{S: S \in (\mathcal{L} \setminus \mathcal{R}) \cup (\mathcal{R} \setminus \mathcal{L})\}$.
Then $m_{\elex T}(\mathcal{L})$ must be strictly larger than
$m_{\elex T}( \mathcal{R})$
since $\{ S \in \mathcal{L}: S \revlex T\}=\{S \in \mathcal{R}:S \revlex T\}$.
Thus we have $T \in \mathcal {L}$.
This fact together with Corollary \ref{conjecture} implies 
that $\Delta(\sigma * \tau) \leq_{Rev} \Delta (\Delta(\sigma)* \Delta(\tau)))$ for all simplicial complexes $\sigma$
and $\tau$,
and Nevo's conjecture is the special case that
$\tau$ consists of two points.

To prove Theorem \ref{bound},
we need some techniques of generic initial ideals
which have a close connection with algebraic shifting.
Let $K$ be an infinite field,
$V$ a $K$-vector space with basis $e_1,e_2,\dots,e_n$
and $E=\bigoplus_{d=0}^n \bigwedge^d V$ the exterior algebra of $V$.
For a graded ideal $J\subset E$,
we write $\Gin_\prec (J)$ for the generic initial ideal of $J$ with respect to a term order $\prec$.
For every monomial $e_S=e_{s_1}\wedge e_{s_2} \wedge \cdots \wedge e_{s_d} \in E$
and for every term order $\prec$,
we write
$$m_{\succeq e_S}(J) = |\{e_R \in J : \ |R|=|S| \ \mbox{and} \ e_R \succeq e_S\}|.$$
We will use the following proposition to prove Theorem \ref{bound}.

\begin{prop}
Let $J\subset E$ be a graded ideal and
$\prec$ and $\prec'$ term orders.
Then, for any monomial $e_S \in E$,
one has
$$m_{\succeq e_S}(\Gin_\prec (J))\geq m_{\succeq e_S}(\Gin_\prec (\init_{\prec'} (J))).$$
\end{prop}

Note that 
the same property as Proposition \ref{prop}
for generic initial ideals over the polynomial ring
was proved by Conca \cite{C}.

This paper is organized as follows:
In \S 1,
we will give the definition of generic initial ideals and recall some basic properties.
In \S 2, we will prove   Proposition \ref{prop}.
In \S 3, we will prove  Theorem \ref{bound} and Corollary \ref{conjecture}.

\section{Generic initial ideals in the exterior algebra}

Let $K$ be an infinite field,
$V$ a $K$-vector space with basis $e_1,e_2,\dots,e_n$
and $E=\bigoplus_{d=0}^n \bigwedge^d V$ the exterior algebra of $V$.
For an integer $d \geq 0$,
let ${[n] \choose d}$ denote the family of $d$-subsets of $[n]$.
If $S= \{ s_1,s_2,\dots,s_d \} \in { [n] \choose d}$ with $s_1<s_2<\cdots<s_d$,
then the element $e_S=e_{s_1} \wedge e_{s_2}\wedge \cdots \wedge e_{s_d}$ will be called a \textit{monomial} of $E$ of degree $d$.
We refer the reader to \cite{AHH} for foundations of the \grev \  theory over the exterior algebra.
Let $\prec$ be a term order.
In this paper,
for $f=\sum_{S \subset [n]} \alpha_S e_S \in E$ with each $\alpha_S \in K$,
we write $\init_\prec (f) = \max_\prec \{ e_S : \alpha_S\ne 0\}$.

Let \GL \ denote the general linear group with coefficients in $K$.
For $\varphi=(a_{ij}) \in \GL$ and for  $f(e_1,\dots,e_n)\in E$,
we define 
$$\varphi (f(e_1,\dots,e_n))=f(\sum_{i=1}^n a_{i1} e_i,\dots,\sum_{i=1}^n a_{in} e_i).$$ 
Also, for a graded ideal $J\subset E$
and for $\varphi \in \GL$, define $\varphi (J)=\{\varphi (f) : f \in J\}$.
A fundamental theorem of generic initial ideals is the following.

\begin{theorem}[{\cite[Theorem 1.6]{AHH}}] \label{gin}
Fix a term order $\prec$. 
Then, for each graded ideal $J\subset E$,
there exists a nonempty Zariski open subset $U\subset \GL$ such that
$\init_\prec( \varphi (J))$ is constant for all $\varphi\in U$.
\end{theorem}

This monomial ideal $\init_\prec(\varphi (J))$ with $\varphi \in U$ is called the \textit{generic initial ideal of} $J$
\textit{with respect to the term order} $\prec$,
and will be denoted $\Gin_\prec(J)$.

\begin{definition} \label{de}
Fix a term order $\prec$.
Given an arbitrary graded ideal $J= \bigoplus_{d=0}^{n} J_d$
of $E$ with each $J_d \subset \bigwedge^dV$,
fix $\varphi \in \GL$ for which
$\init_{\prec}(\varphi(J))$
is the generic initial ideal $\Gin_\prec(J)$ of $J$.
Recall that the subspace $\bigwedge^dV$
is of dimension ${n \choose d}$
with a canonical $K$-basis
$\{e_S:S \in {[n] \choose d}\}$.
Choose an arbitrary $K$-basis $f_1, \ldots, f_m$ of
$J_d$, where $m = \dim_K J_d$.  Write each $\varphi(f_i)$,
$1 \leq i \leq m$, of the form
\[
\varphi(f_i) = \sum_{S \in {[n] \choose d}}
\alpha_{i}^{S} \, e_{S}
\]
with each $\alpha_{i}^{S} \in K$.
Let $M(J,d)$ denote the $m \times {n \choose d}$ matrix
\[
M(J,d) =(\alpha_{i}^{S})_{1 \leq i \leq m, \,
S \in {[n] \choose d}}
\]
whose columns are indexed by
$S \in {[n] \choose d}$.
For each $S \in {[n] \choose d}$,
write
$\sub$ for the submatrix of $M(J,d)$
which consists of the columns of $M(J,d)$ indexed by
those $R \in {[n] \choose d}$ with
$R \succeq S$
and write
$\tsub$ for the submatrix of $\sub$
which is obtained by removing the column
of $\sub$ indexed by $S$.

It is not hard to see that
we can know generic initial ideals by using the rank of these matrices.
Indeed, following properties are known.
\end{definition}

\begin{lemma}[{\cite[Lemma 2.1]{MH}}] \label{generic}
Let $e_{S} \in \bigwedge^dV$
with $S \in {[n] \choose d}$.
Then
one has $e_{S} \in (\Gin_\prec(J))_d$
if and only if
$\rank(\tsub) < \rank(\sub)$.
\end{lemma}

\begin{lemma}[{\cite[Corollary 2.2]{MH}}] \label{independent}
The rank of a matrix $\sub$,
$S \in {[n] \choose d}$,
is independent of
the choice of $\varphi \in \GL$
for which
$\Gin_\prec(J) = \init_{\prec}(\varphi(J))$
and independent of the choice of the $K$-basis $f_1, \ldots, f_m$ of $J_d$.
\end{lemma}

\begin{lemma}[{\cite[Corollary 2.3]{MH}}]
\label{isomorphic}
Let $J \subset E$ be a graded ideal
and $\psi \in \GL$.
Then one has
$\rank(\sub) = \rank(M_{\succeq S}(\psi(J),d))$
for all $S \in {[n] \choose d}$.
\end{lemma}

Also, the next lemma immediately follows from Lemma \ref{generic}.

\begin{lemma}
\label{rank}
Let $J\subset E$ be a graded ideal.
For every $S\in { [n] \choose d}$,
one has
\[
m_{\succeq e_S}(\Gin_\prec(J)) = \rank (\sub).
\]
\end{lemma}

\section{proof of proposition \ref{prop}}

We will follow the basic technique developed in \cite{G}.
(See also \cite[Chapter 15]{E}.)


\begin{lemma}[{\cite[Corollary 1.7]{G}}] \label{order}
For any term order $\prec$
and for any finite set of monomials $M \subset E$,
there exist positive integers $d_1,d_2,\dots ,d_n$ such that 
for any $e_{S},e_{R}\in M$ with $|S|=|R|$, one has
$e_{S}\prec e_{R} $ if and only if $\sum_{k\in S} d_k > \sum_{k\in R} d_k$.
\end{lemma}

For every ideal $J\subset E$,
a subset $G=\{g_1,\dots,g_m\}\subset J$ is called a \textit{\grev \ of} $J$ with respect to $\prec$
if  $\{\init_\prec(g_1),\dots,\init_\prec(g_m)\}$ generates $\init_\prec(J)$.
A \grev \  always exists and is actually a generating set of $J$
(\cite[Theorem 1.4]{AHH}).

\begin{lemma}\label{flat}
Let $K[t]$ be the polynomial ring.
Fix a term order $\prec$.
For every graded ideal $J\subset E$,
there is a subset $G(t)=\{g_1(t),\dots,g_m(t)\} \subset  E\otimes K[t]$ which satisfies
the following conditions.
\begin{itemize}
\item[(i)] One has
$  g_i(0)=\init_\prec (g_i(t_0))$
 for all $t_0 \in K$;
\item[(ii)] Let  $J(t_0)$ with $t_0 \in K$ be the ideal generated by $G(t_0)$.
If $t_0 \ne 0$, then there exists  $\varphi_{t_0}\in \GL$
 such that $\varphi_{t_0}(J)=J(t_0)$;
\item[(iii)] For all $t_0 \in K$, $G(t_0)$ is a \grev \ of $J(t_0)$ with respect to $\prec$ and
$\init_\prec(J(t_0)) = \init_\prec(J)$.  
\end{itemize}
\end{lemma} 
\begin{proof}
Let $G=\{g_1,g_2,\dots,g_m \}$ be a \grev \ of $J$ with respect to $\prec$,
where each $g_j$ is homogeneous.
Let $M\subset E$ be the set of monomials.
Since $M$ is a finite set,
Lemma \ref{order} says that there exist positive integers $d_1,d_2,\dots,d_n$ such that,
for any $e_{S},e_{R} \in M$ with $|S|=|R|$,
\begin{eqnarray}
e_{S} \prec e_{R} \mbox{ if and only  if } \sum_{k\in S} d_k  > \sum_{k \in R} d_k. \label{maru1}
\end{eqnarray}
Let $e_{S_i} = \init_\prec (g_i)$.
For each $R \subset [n]$, write $d(R)=  \sum_{k\in R} d_k$.
Set
$$g_i(e_1,\dots,e_n)(t) = t^{-d( S_i)} g_i (t^{d_1}e_1,\dots,t^{d_n}e_n)$$
and $G(t)=\{  g_1(t),\dots, g_m(t)\}$.
We will show that this set $G(t)$ satisfies conditions (i), (ii) and (iii).

First, we will show (i).
Each $g_i(t)$ can be written in the form
$$g_i(t) = \alpha_{S_i} e_{S_i} + \sum_{R \prec S_i} \alpha_R \cdot t^{ \{ d(R)-d(S_i)\}} \cdot e_R,$$
where $\alpha_{S_i} \in K \setminus \{0\}$ and each $\alpha_R \in K$.
Then (\ref{maru1}) says that $d(R)-d(S_i)>0$ for all $R$ with $\alpha_R \ne 0$.
Thus we have $ g_i(0)=\init_\prec(g_i) = \init_\prec ( g_i(t_0))$ for all $t_0 \in K$ as desired.

Second,
for each $t_0\in K \setminus \{0\}$,
define a matrix $\varphi_{t_0} \in \GL$ by 
$$\varphi_{t_0}(e_i)={t_0}^{d_i}e_i \ \ \ \ \ \ \ \mbox{ for } i=1,2,\dots,n.$$
Then 
the construction of $g_i(t)$ says that
$\varphi_{t_0}(g_i)=  t_0^{d(S_i)}  g_i(t_0)$,
and therefore we have
$\varphi_{t_0} (J) =J(t_0)$ for all $t_0\in K \setminus \{0\}$.
Thus (ii) is satisfied.

Finally, we will show (iii).
Since we already proved $\init_\prec (g_i(t_0))= \init_\prec(g_i)$ for all $t_0 \in K$,
what we must prove is $\init_\prec ( J(t_0)) = \init_\prec (J)$.
The inclusion $\init_\prec ( J(t_0)) \supset \init_\prec (J)$ follows from $\init_\prec (g_i(t_0))= \init_\prec(g_i)$.
Recall that $J$ and $\init_\prec(J)$ have the same Hilbert function,
i.e.,
we have $\dim_K(J_d)=\dim_K(\init_\prec(J)_d)$ for all $d>0$.
Then we have
$$\dim_K (\init_\prec(J(t_0))_d)=\dim_K(J(t_0)_d)=\dim_K(J_d)=\dim_K(\init_\prec(J)_d).$$
Hence we have $\init_\prec ( J(t_0)) = \init_\prec (J)$.
Thus $G(t_0)$ is a \grev \ of $J(t_0)$  for all $t_0\in K$.
\end{proof}

\begin{lemma}\label{flatdegree}
Let $J\subset E$ be a graded ideal.
For every $t_0\in K$,
let $J(t_0) \subset  E$ be the ideal given in Lemma \ref{flat}.
Then, for all $d>0$, there exists a subset $G_d(t)\subset E\otimes K[t]$ such that
$G_d(t_0)$ is a $K$-basis of $J(t_0)_d$ for all $t_0\in K$.
\end{lemma}

\begin{proof}
Let $G(t)=\{g_1(t),\dots,g_m(t)\}$ be a subset of $E\otimes K[t]$  which is given in Lemma \ref{flat}.
Let
$$\widetilde G_d(t)=\{e_S \wedge g_i(t) : \deg(e_S \wedge g_i(0)) =d \mbox{ and } e_S \wedge g_i(0) \ne 0  \}.$$
For every $t_0 \in K$,
since $G(t_0)$ is a \grev \ of $J(t_0)$
and $g_i(0)=\init_\prec(g_i(t_0))$,
the set $\{ \init_\prec (h(t_0)): h(t_0) \in \widetilde G_d(t_0) \}$ spans $\init_\prec (J(t_0))_d=\init_\prec (J)_d$.
Also, Lemma \ref{flat} (i) says that $h(0) = \init_\prec (h(t_0))$ for all $t_0 \in K$.
Thus there is a subset $G_d(t)\subset \widetilde G_d(t)$ such that
$G_d(0)$ is a $K$-basis of $\init_\prec (J)_d$.
On the other hand, for any $t_0 \in K$,
since each $h(t_0) \in G_d(t_0)$ has a different initial monomial,
the set $G_d(t_0)$ is linearly independent.
Thus we have 
$$\dim_K(\mathrm{span} \{G_d(t_0)\})=\dim_K(\mathrm{span} \{G_d(0)\})=\dim_K(\init_\prec (J)_d)=\dim_K( J(t_0)_d),$$
where $\mathrm{span}(A)$ denotes the $K$-vector space spanned by a finite set $A \subset E$.
Hence  $G_d(t_0)$ is a $K$-basis of $J(t_0)_d$ for all $t_0\in K$.
\end{proof}

\begin{proposition} \label{prop}
Let $J\subset E$ be a graded ideal and
$\prec$ and $\prec'$ term orders.
Then, for any monomial $e_S \in E$,
one has
$$m_{\succeq e_S}(\Gin_\prec (J))\geq m_{\succeq e_S}(\Gin_\prec (\init_{\prec'} (J))).$$
\end{proposition}
\begin{proof}
First,
by Lemma \ref{rank}, we have 
$$m_{\succeq e_S}(\Gin_\prec (J)) = \rank( \sub)$$
and
$$m_{\succeq e_S}(\Gin_\prec (\init_{\prec'} (J))) = \rank( M_{\succeq S}(\init_{\prec'}(J),d)).$$
Thus what we must prove is
$\rank( \sub) \geq \rank  (M_{\succeq S} (\init_{\prec'} (J),d ))$.

Let $m=\dim_K (J_d)$,
$G_d(t)=\{g_1(t),\dots,g_m(t)\} \subset  E\otimes K[t]$ the subset given in 
Lemma \ref{flatdegree} w.r.t.\ the term order $\prec'$
and $J(t_0)$, where $t_0 \in K$, the ideal given in Lemma \ref{flat}.
Then, for each $t_0 \in K \setminus \{0\}$,
there exists $\varphi_{t_0}\in\GL$ such that $\varphi_{t_0}(J) = J(t_0)$.
Thus Lemma \ref{isomorphic} says that we have
\begin{eqnarray}
\mathrm{rank} (\sub )= \mathrm{rank} (M_{\succeq S} (J(t_0),d)) \ \ \mbox{  for all } t_0 \in K \setminus \{0\}. \label{iti}
\end{eqnarray}

Let $A\subset K$ be a finite set with $0\in A$  and $|A| \gg 0$.
Then Theorem \ref{gin} says that, for each $a \in A$, 
there exists a nonempty Zariski open subset $U_a \subset \GL$
such that $\Gin (J(a))= \init_\prec (\varphi (J(a)))$ for all $\varphi \in U_a$.
Since $U= \bigcap_{ a \in A} U_a$ is also a nonempty Zariski open subset of $\GL$,
we have
$\Gin_\prec(J(a))=\init_\prec (\varphi (J(a)))$,
for all $\varphi \in U$ and all $a \in A$. 

Fix $\varphi \in U$.
Each $\varphi(g_i(t))$, where $1\leq i \leq m$, can be written in the form
$$\varphi(g_i(t))=\sum_{S\in {[n] \choose d} }\alpha_i^{S}(t) e_S,$$
where $\alpha_i^{S}(t) \in K[t]$.
Define the matrix $\lsub t=(\alpha_i^R(t))_{1\leq i\leq m,\ R\succeq S}$
by the same way as Definition \ref{de}.
Recall that Lemma \ref{flatdegree} says that 
$G_d(a)$ is a $K$-basis of $J(a)_d$ for all $a \in A$.
Since Lemma \ref{independent} says that  $\rank (M_{\succeq S} (J(a),d))$ is independent of the choice of a $K$-basis
and independent of the choice of $\varphi \in U_a$,
it follows that
\begin{eqnarray}
\mathrm{rank}(\lsub a)=\mathrm{rank}(M_{\succeq S} (J(a),d))\ \ \ \mbox{ for all } a \in A. \label{ni}
\end{eqnarray}
Let $l=\max \{\deg(\alpha_i^S (t) ): 1\leq i \leq m, \ S \in {[n] \choose d} \}$.
Recall that the rank of matrices is equal to the maximal size of nonzero minors.
In this case, each minor of $\lsub t$ is a polynomial of $K[t]$
and has at most degree $l m$.
Furthermore, the number of nonzero minors of $\lsub t$ is finite.
Since integers $l$ and $m$ do not depend on $A$ and
$|A|$ is sufficiently large,
there exists $a_0\ne 0 \in A$ such that $f(a_0) \ne 0$ for all nonzero minors $f(t)$ of $\lsub t$.
In particular, we have
$\rank (\lsub {a_0})\geq \rank (\lsub {t_0})$ for all $t_0\in K$.
Recall that Lemma \ref{flat} (i) and (iii) say that $J(0)= \init_{\prec'}(J)$.
Then, by (\ref{iti}) and (\ref{ni}), we have 
\begin{eqnarray*}
\rank( \sub)&=&
\rank( M_{\succeq S} (J(a_0),d))\\
&=&\rank (\lsub {a_0})\\
&\geq&  \rank ( \lsub 0) = \rank  (M_{\succeq S} (\init_{\prec'} (J),d )),
\end{eqnarray*}
as desired.
\end{proof}

\section{exterior shifting of the join of simplicial complex}
Let $\sigma$ be a simplicial complex on $[n]$.
The \textit{exterior face ideal $J_\sigma$ of $\sigma$} is the ideal of $E$ generated by all monomials $e_S$ with $S \not\in \sigma$.
%
For every $\varphi \in \GL$ and for every simplicial complex $\sigma$,
the simplicial complex $\Delta_{\varphi} (\sigma)$ is defined by
 $J_{\Delta_\varphi (\sigma) }=\init_\revlex (\varphi (J_\sigma))$.
The \textit{exterior algebraic shifted complex $\Delta(\sigma)$ of $\sigma$} is the simplicial complex defined by
$J_{\Delta(\sigma)} = \Gin_{\revlex} (J_\sigma)$.

\begin{theorem} \label{bound}
Let $\sigma$ be a simplicial complex on $[n]$ and
$\varphi \in \GL$.
Then, for any $S\subset [n]$,
one has
$$m_{\elex S}(\Delta(\sigma)) \geq m_{\elex S}(\Delta (\Delta_{\varphi}(\sigma))).$$
\end{theorem}

\begin{proof}
By the definition of exterior face ideals,
for any $d$-subset $S\in {[n] \choose d}$ and for any simplicial complex $\tau$,
we have 
\begin{eqnarray}
m_{\elex S}(\tau)&=& |\tau_{d-1}|-|\{R\in \tau_{d-1} : R \succ_{rev} S\}| \nonumber \\
&=&|\tau_{d-1} |- |\{ R \in { [n] \choose d}: R \succ_{rev} S\}|  + 
|\{ e_R\in I _\tau: R\succ_{rev} S\}|.
\label{hii}
\end{eqnarray}
On the other hand,  Proposition \ref{prop} says that, for any $S \subset [n]$, one has
\begin{eqnarray}
m_{\rlex e_S}(\Gin_\revlex (\varphi(J_\sigma)))\geq m_{\rlex e_S}(\Gin_\revlex ( \init_\revlex (\varphi(J_\sigma)))). \label{huu}
\end{eqnarray}
Also, Lemma \ref{isomorphic} says that
\begin{eqnarray}
J_{\Delta(\sigma)}=\Gin_{\revlex} (J_\sigma)=\Gin_{\revlex}(\varphi(J_\sigma))  \label{mii}
\end{eqnarray}
and
\begin{eqnarray}
 J_{\Delta(\Delta_{\varphi}(\sigma))}=\Gin_{\revlex}(\init_{\revlex} (\varphi(J_\sigma))). \label{you}
\end{eqnarray}
Then, equalities (\ref{hii}), (\ref{huu}), (\ref{mii}) and (\ref{you}) say that
$$m_{\elex S}(\Delta(\sigma)) \geq m_{\elex S}(\Delta (\Delta_\varphi(\sigma))),$$
as desired.
\end{proof}

\begin{cor}\label{conjecture}
Let $\sigma$ be a simplicial complex on $\{1,2,\dots,k\}$ and
$\tau$ a simplicial complex on $\{k+1,k+2,\dots,n\}$.
Then, for any $S\subset [n]$, one has 
$$m_{\elex S} (\Delta(\sigma * \tau))\geq m_{\elex S} (\Delta(\Delta(\sigma) * \Delta(\tau))).$$
\end{cor}

\begin{proof}
Let $l=|\{k+1,k+2,\dots,n\}|.$
Then there exist $\varphi\in\mathrm{GL}_k(K)$ and $\psi \in \mathrm{GL}_l(K)$ 
such that
$\Delta_\varphi(\sigma)=\Delta(\sigma)$ and $\Delta_\psi(\tau)=\Delta(\tau)$.
For any $\varphi\in\mathrm{GL}_k(K)$, we define $\bar \varphi \in \GL$ by
\begin{eqnarray*}
\bar \varphi (e_i) =
\left\{
\begin{array}{l}
\varphi(e_i), \ \ \mbox{if} \ i\in \{1,2,\dots,k\}, \\
e_i, \hspace{24pt} \mbox{otherwise}.
\end{array}
\right.
\end{eqnarray*}
Also, for any $\psi \in \mathrm{GL}_l(K)$, define $\bar \psi \in \GL$ in the same way.
Then we have
$$\Delta_{\bar \varphi}(\Delta_{\bar \psi} (\sigma * \tau))
=\Delta_{\bar \varphi}(\sigma *\Delta(\tau))
=\Delta(\sigma) * \Delta(\tau).$$
Then Theorem \ref{bound} says that
$$m_{\elex S} (\Delta(\sigma * \tau))\geq
 m_{\elex S} (\Delta(\Delta_{\bar \varphi}(\Delta_{\bar \psi} (\sigma * \tau))))=
 m_{\elex S} (\Delta(\Delta(\sigma) * \Delta(\tau))),$$
as desired.
\end{proof}

\begin{rem}
Proposition \ref{prop} holds for an arbitrary term order.
However, Theorem \ref{bound} and Corollary \ref{conjecture} only hold for the reverse lexicographic order.
Recall that Nevo's example says that
$$\Delta(\Sigma(\sigma))= \{\{1,2,3\},\{1,2,4\},\{1,2,5\},\{1,2,6\}\}$$
and
$$\Delta(\Sigma(\Delta(\sigma)))= \{\{1,2,3\},\{1,2,4\},\{1,2,5\},\{1,3,4\}\},$$
where $\sigma$ is the simplicial complex generated by $\{1,2\}$ and $\{3,4\}$.
In this case, 
$\{1,2,5\}$ and $\{1,2,6\}$ are larger than $\{1,3,4\}$ with respect to the lexicographic order
$\prec_{lex}$ induced by $1<2<\cdots<n$.
This implies that $m_{ \preceq_{lex} \{1,3,4\}} ( \Delta(\Sigma(\sigma)))=2$
but  $m_{ \preceq_{lex} \{1,3,4\}} (\Delta(\Sigma(\Delta(\sigma))))=3$.
\end{rem}


\begin{thebibliography}{1}
\bibitem{AHH}
A. Aramova, J. Herzog and T. Hibi,
Gotzmann theorems for exterior algebras and combinatorics,
 {\em J. Algebra} \textbf{191} (1997), 174--211.


\bibitem{BK}
A. B\"ojrner  and G. Kalai,
An extended Euler-Poincar\'e theorem,
\textit{Acta Math.}
\textbf{161} (1988), 279--303.


\bibitem{C}
A. Conca,
Reduction numbers and initial ideals,
\textit{Proc. Amer. Math. Soc.}
\textbf{131} (2003), 1015--1020.

\bibitem{E} 
D. Eisenbud,
\textit{Commutative algebra, with view towards algebraic geometry,}
Graduate Texts in Mathematics \textbf{150},
Spriger-Verlag, New York, 1995.




\bibitem{G} M. Green, 
Generic initial ideals, \textit{in} ``Six Lectures on Commutative Algebra"
(J. Elias, J. M. Giral, R. M. Mir'o-Roig, and S. Zarzuela, Eds.),
Progress in Math., \textbf{166}, Birkh\"auser, Basel, 1998, 119--186.

\bibitem{H} J. Herzog, Generic initial ideals and graded Betti numbers,
{\em in} ``Computational Commutative Algebra and Combinatorics''
(T. Hibi, Ed.), Advanced Studies in Pure Math., Volume 33, 2002,
pp. 75--120.


\bibitem{K}
G. Kalai, Algebraic shifting,
{\em in} ``Computational Commutative Algebra and Combinatorics''
(T. Hibi, Ed.), Advanced Studies in Pure Math., Volume 33, 2002,
pp. 121--163.


\bibitem{MH}
S. Murai and T. Hibi,
The behavior of graded Betti numbers via algebraic shifting and combinatorial shifting,
Preprint, 2005, {arXiv:math.AC/0503685}.


\bibitem{N} E. Nevo,
Algebraic shifting and basic constructions on simplicial complexes,
\textit{J. Algebraic Combin.} \textbf{22} (2005), 411--433.

\end{thebibliography}
\end{document}